\documentclass[12pt,oneside,french,english]{amsart}
\usepackage[T1]{fontenc}
\usepackage[latin9]{inputenc}
\usepackage[letterpaper]{geometry}
\usepackage{amsthm}
\usepackage{amssymb}
\usepackage{esint}

\makeatletter

\newcommand{\noun}[1]{\textsc{#1}}

\numberwithin{equation}{section}
\numberwithin{figure}{section}

\makeatother

\usepackage{babel}
\addto\extrasfrench{\providecommand{\fg}{\ifdim\lastskip>\z@\unskip\fi~\frqq}}

\begin{document}
\selectlanguage{french}%
\textbf{\hfill{}} \textit{To appear in Journal of Mathematical Sciences}

\selectlanguage{english}%
\vspace{0.5cm}

\title{Application of a Bernstein type inequality to rational interpolation
in the Dirichlet space}

\author{Rachid Zarouf}
\begin{abstract}
We prove a Bernstein-type inequality involving the Bergman and the
Hardy norms, for rational functions in the unit disc $\mathbb{D}$
having at most $n$ poles all outside of $\frac{1}{r}\mathbb{D}$,
$0<r<1$. The asymptotic sharpness of this inequality is shown as
$n\rightarrow\infty$ and $r\rightarrow1^{-}.$ We apply our Bernstein-type
inequality to an effective Nevanlinna-Pick interpolation problem in
the standard Dirichlet space, constrained by the $H^{2}$- norm.
\end{abstract}
\maketitle

\section*{Introduction}

\subsection*{a. Statement of the problems}

\begin{flushleft}
Let $\mathbb{D}=\{z\in\mathbb{C}:\,\vert z\vert<1\}$ be the unit
disc of the complex plane and let ${\rm Hol}\left(\mathbb{D}\right)$
be the space of holomorphic functions on $\mathbb{D}.$ Let also $X$
and $Y$ be two Banach spaces of holomorphic functions on the unit
disc $\mathbb{D},$ $X,\, Y\subset{\rm Hol}\left(\mathbb{D}\right).$
Here and later on, $H^{\infty}$ stands for the space (algebra) of
bounded holomorphic functions in the unit disc $\mathbb{D}$ endowed
with the norm $\left\Vert f\right\Vert _{\infty}=\sup_{z\in\mathbb{D}}\left|f(z)\right|.$
We suppose that $n\geq1$ is an integer, $r\in[0,\,1)$ and we consider
the two following problems.
\par\end{flushleft}

\vspace{0.1cm}

\textbf{Problem 1.} Let $\mathcal{P}_{n}$ be the complex space of
analytic polynomials of degree less or equal than $n$, and \[
\mathcal{R}_{n,\, r}=\left\{ \frac{p}{q}\,:\; q\in\mathcal{P}_{n},\; d^{\circ}p<d^{\circ}q,\; q(\zeta)=0\Longrightarrow\left|\zeta\right|\geq\frac{1}{r}\right\} ,\]
(where $d^{\circ}p$ means the degree of any $p\in\mathcal{P}_{n}$)
be the set of all rational functions in $\mathbb{D}$ of degree less
or equal than $n\ge1$, having at most $n$ poles all outside of $\frac{1}{r}\mathbb{D}.$
Notice that for $r=0$, we get $\mathcal{R}_{n,\,0}=\mathcal{P}_{n-1}$.
Our first problem is to search for the {}``best possible'' constant
$\mathcal{C}_{n,\, r}(X,\, Y)$ such that \[
\left\Vert f'\right\Vert _{X}\leq\mathcal{C}_{n,\, r}(X,\, Y)\left\Vert f\right\Vert _{Y}\]
for all $f\in\mathcal{R}_{n,\, r}.$ 

\vspace{0.1cm}

\textbf{Problem 2.} Let $\sigma=\left\{ \lambda_{1},...,\lambda_{n}\right\} $
be a finite subset of $\mathbb{D}$. What is the best possible interpolation
by functions of the space $Y$ for the traces $f_{\vert\sigma}$ of
functions of the space $X$, in the worst case? The case $X\subset Y$
is of no interest, and so one can suppose that either $Y\subset X$
or $X$ and $Y$ are incomparable. More precisely, our second problem
is to compute or estimate the following interpolation constant\[
I\left(\sigma,\, X,\, Y\right)=\sup_{f\in X,\,\parallel f\parallel_{X}\leq1}\mbox{inf}\left\{ \left\Vert g\right\Vert _{Y}:\, g_{\vert\sigma}=f_{\vert\sigma}\right\} .\]
We also define \[
\mathcal{I}_{n,\, r}(X,\, Y)=\mbox{sup}\left\{ I(\sigma,\, X,\, Y)\,:\,{\rm card}\,\sigma\leq n\,,\,\left|\lambda\right|\leq r,\,\forall\lambda\in\sigma\right\} .\]

\subsection*{b. Motivations}

$\,$

\textbf{Problem 1. }Bernstein-type inequalities for rational functions
are applied

\textbf{1.1.} in matrix analysis and in operator theory (see {}``Kreiss
Matrix Theorem'' {[}LeTr, Sp{]} or {[}Z1, Z4{]} for resolvent estimates
of power bounded matrices),

\textbf{1.2.} to {}``inverse theorems of rational approximation''
using the \textit{classical Bernstein decomposition} (see {[}Da, Pel,
Pek{]}), 

\textbf{1.3.} to effective $H^{\infty}$ interpolation problems (see
{[}Z3{]} and our Theorem B below in Subsection d), and more generally
to our Problem 1.

\vspace{0.1cm}

\textbf{Problem 2.} We can give three main motivations for Problem
2. 

\textbf{2.1.} It is explained in {[}Z3{]} (the case $Y=H^{\infty})$
why the classical interpolation problems, those of Nevanlinna-Pick
(1908) and Carathéodory-Schur (1916) (see {[}N2{]} p.231 for these
two problems), on the one hand and Carleson's free interpolation problem
(1958) (see {[}N1{]} p.158) on the other hand, are of the nature of
our interpolation problem.

\textbf{2.2.} It is also explained in {[}Z3{]} why this constrained
interpolation is motivated by some applications in matrix analysis
and in operator theory.

\textbf{2.3.} It has already been proved in {[}Z3{]} that for $X=H^{2}$
(see Subsection c. for the definition of $H^{2}$) and $Y=H^{\infty},$ 

\def\theequation{${1}$}\begin{equation}
\frac{1}{4\sqrt{2}}\frac{\sqrt{n}}{\sqrt{1-r}}\leq\mathcal{I}_{n,\, r}\left(H^{2},\, H^{\infty}\right)\leq\sqrt{2}\frac{\sqrt{n}}{\sqrt{1-r}}.\label{eq:-2-1}\end{equation}
The above estimate (1) answers a question of L. Baratchart (private
communication), which is part of a more complicated question arising
in an applied situation in {[}BL1{]} and {[}BL2{]}: given a set $\sigma\subset\mathbb{D}$,
how to estimate $I\left(\sigma,\, H^{2},\, H^{\infty}\right)$ in
terms of $n=\mbox{card}(\sigma)$ and $\mbox{max}{}_{\lambda\in\sigma}\left|\lambda\right|=r$
only?

\subsection*{c. The spaces $X$ and $Y$ considered here}

\begin{flushleft}
Now let us define some Banach spaces $X$ and $Y$ of holomorphic
functions in $\mathbb{D}$ which we will consider throughout this
paper. From now on, if $f\in{\rm Hol}(\mathbb{D})$ and $k\in\mathbb{N}$,
$\hat{f}(k)$ stands for the $k^{th}$ Taylor coefficient of $f.$
\par\end{flushleft}

\textbf{1.} The standard Hardy space $H^{2}=H^{2}(\mathbb{D}),$ \[
H^{2}=\left\{ f\in{\rm Hol}\left(\mathbb{D}\right):\:\left\Vert f\right\Vert _{H^{2}}^{2}=\sup_{0\leq r<1}\int_{\mathbb{T}}\left|f(rz)\right|^{2}{\rm d}m(z)<\infty\right\} ,\]
 where $m$ stands for the normalized Lebesgue measure on $\mathbb{T}=\{z\in\mathbb{C}:\,\vert z\vert=1\}.$
An equivalent description of the space $H^{2}$ is \[
H^{2}=\left\{ f=\sum_{k\geq0}\hat{f}(k)z^{k}:\,\,\left\Vert f\right\Vert _{H^{2}}=\left(\sum_{k\geq0}\left|\hat{f}(k)\right|^{2}\right)^{\frac{1}{2}}<\infty\right\} .\]

\textbf{2}. The standard Bergman space $L_{a}^{2}=L_{a}^{2}\left(\mathbb{D}\right),$
\[
L_{a}^{2}=\left\{ f\in{\rm Hol}\left(\mathbb{D}\right):\:\left\Vert f\right\Vert _{L_{a}^{2}}^{2}=\frac{1}{\pi}\int_{\mathbb{D}}\left|f(z)\right|^{2}{\rm d}A(z)<\infty\right\} ,\]
where $A$ is the standard area measure, also defined by \[
L_{a}^{2}=\left\{ f=\sum_{k\geq0}\hat{f}(k)z^{k}:\,\left\Vert f\right\Vert _{L_{a}^{2}}\,=\left(\sum_{k\geq0}\left|\hat{f}(k)\right|^{2}\frac{1}{k+1}\right)^{\frac{1}{2}}<\infty\right\} .\]

\textbf{3.} The analytic Besov space $B_{2,\,2}^{\frac{1}{2}}$ (also
known as the standard Dirichlet space) defined by \[
B_{2,\,2}^{\frac{1}{2}}=\left\{ f=\sum_{k\geq0}\hat{f}(k)z^{k}:\,\left\Vert f\right\Vert _{B_{2,\,2}^{\frac{1}{2}}}=\left(\sum_{k\geq0}(k+1)\left|\hat{f}(k)\right|^{2}\right)^{\frac{1}{2}}<\infty\right\} .\]
 Then if $f\in B_{2,\,2}^{\frac{1}{2}},$ we have the following equality

\def\theequation{${2}$}\begin{equation}
\left\Vert f\right\Vert _{B_{2,\,2}^{\frac{1}{2}}}^{2}=\left\Vert f'\right\Vert _{L_{a}^{2}}^{2}+\left\Vert f\right\Vert _{H^{2}}^{2},\label{eq:-2}\end{equation}
which establishes a link between the spaces $B_{2,\,2}^{\frac{1}{2}}$
and $L_{a}^{2}$.

\subsection*{d. The results}

Here and later on, the letter $c$ denotes a positive constant that
may change from one step to the next. For two positive functions $a$
and $b$, we say that $a$ is dominated by $b$, denoted by $a=O(b),$
if there is a constant $c>0$ such that $a\leq cb;$ and we say that
$a$ and $b$ are comparable, denoted by $a\asymp b$, if both $a=O(b)$
and $b=O(a)$ hold.

\vspace{0.1cm}

\textbf{Problem 1.} Our first result (Theorem A, below) is a partial
case ($p=q=2$, $s=\frac{1}{2}$) of the following K. Dyakonov's result
{[}Dy{]}: if $p\in[1,\,\infty),$ $s\in(0,\,+\infty),\; q\in[1,\,+\infty]$,
then there exists a constant $c_{p,\, s}>0$ such that

\def\theequation{${3}$}\begin{equation}
\mathcal{C}_{n,\, r}\left(B_{p,\, p}^{s-1},\, H^{q}\right)\leq c_{p,\, s}\sup\left\Vert B'\right\Vert _{H^{\gamma}}^{s},\label{eq:-5}\end{equation}
where $\gamma$ is such that $\frac{s}{\gamma}+\frac{1}{q}=\frac{1}{p},$
and the supremum is taken over all finite Blaschke products $B$ of
order $n$ with $n$ zeros outside of $\frac{1}{r}\mathbb{D}.$ Here
$B_{p,\, p}^{s}$ stands for the Hardy-Besov space which consists
of analytic functions $f$ on $\mathbb{D}$ satisfying\[
\left\Vert f\right\Vert _{B_{p,\, p}^{s}}=\sum_{k=0}^{n-1}\left|f^{(k)}(0)\right|+\int_{\mathbb{D}}\left(1-\left|w\right|\right)^{(n-s)p-1}\left|f^{(n)}(w)\right|^{p}{\rm d}A(w)<\infty.\]
 For the (tiny) partial case considered here, our proof is different
and the constant $c_{2,\,\frac{1}{2}}$ is asymptotically sharp as
$r$ tends to $1^{-}$ and $n$ tends to $+\infty$. 

\begin{flushleft}
\textbf{Theorem A.} \textit{Let $n\geq1$ and $r\in[0,\,1).$ We have}
\par\end{flushleft}

\textit{(i)}

\def\theequation{${4}$}

\textit{\begin{equation}
\widetilde{a}(n,\, r)\sqrt{\frac{n}{1-r}}\leq\mathcal{C}_{n,\, r}\left(L_{a}^{2},\, H^{2}\right)\leq\widetilde{A}(n,\, r)\sqrt{\frac{n}{1-r}},\label{eq:-4}\end{equation}
where \[
\widetilde{a}(n,\, r)\geq\left(1-\frac{1-r}{n}\right)^{\frac{1}{2}}\; and\;\widetilde{A}(n,\, r)\leq\left(1+r+\frac{1}{\sqrt{n}}\right)^{\frac{1}{2}}.\]
}

\textit{(ii) Moreover, the sequence \[
\left(\frac{\mathcal{C}_{n,\, r}\left(L_{a}^{2},\, H^{2}\right)}{\sqrt{n}}\right)_{n\geq1}\]
is convergent and there exists a limit}

\textit{\def\theequation{${5}$}\begin{equation}
\lim_{n\rightarrow\infty}\frac{\mathcal{C}_{n,\, r}\left(L_{a}^{2},\, H^{2}\right)}{\sqrt{n}}=\sqrt{\frac{1+r}{1-r}}.\label{eq:-4}\end{equation}
for all $r\in[0,\,1)$. }

$\,$

Notice that it has already been proved in {[}Z2{]} that there exists
a limit

\def\theequation{${6}$}\begin{equation}
\lim_{n\rightarrow\infty}\frac{\mathcal{C}_{n,\, r}\left(H^{2},\, H^{2}\right)}{n}=\frac{1+r}{1-r},\label{eq:-6}\end{equation}
for every $r,\;0\leq r<1$. 

\vspace{0.1cm}

\textbf{Problem 2.} Looking at motivation 2.3, we replace the algebra
$H^{\infty}$ by the Dirichlet space $B_{2,\,2}^{\frac{1}{2}}.$ We
show that the {}``gap'' between $X=H^{2}$ and $Y=H^{\infty}$ (see
(1)) is asymptotically the same as the one which exists between $X=H^{2}$
and $Y=B_{2,\,2}^{\frac{1}{2}}.$ In other words,

\def\theequation{${7}$}\textit{\begin{equation}
\mathcal{I}_{n,\, r}\left(H^{2},\, B_{2,\,2}^{\frac{1}{2}}\right)\asymp\mathcal{I}_{n,\, r}\left(H^{2},\, H^{\infty}\right)\asymp\sqrt{\frac{n}{1-r}}.\label{eq:-1}\end{equation}
}More precisely, we prove the following Theorem B, in which the right-hand
side inequality of (10) is a consequence of the right-hand side inequality
of (4) in the above Theorem A. 

\textbf{Theorem B. }\textit{Let $n\geq1$, and $r\in[0,\,1).$ Then,}

\def\theequation{${8}$}\textit{\begin{equation}
\mathcal{I}_{n,\, r}\left(H^{2},\, B_{2,\,2}^{\frac{1}{2}}\right)\leq\left[\left(\mathcal{C}_{n,\, r}\left(L_{a}^{2},\, H^{2}\right)\right)^{2}+1\right]^{\frac{1}{2}}.\label{eq:-3}\end{equation}
}

\begin{flushleft}
\textit{Let $\lambda\in\mathbb{D}$ and the corresponding one-point
interpolation set $\sigma_{n,\,\lambda}=\underbrace{\{\lambda,\lambda,...,\lambda\}}_{n}.$
We have,}
\par\end{flushleft}

\def\theequation{${9}$}\textit{\begin{equation}
I\left(\sigma_{n,\,\lambda},H^{2},\, B_{2,\,2}^{\frac{1}{2}}\right)\geq\sqrt{\frac{n}{1-\left|\lambda\right|}}\left[\frac{(1+\left|\lambda\right|)^{2}-\frac{2}{n}-\frac{2\left|\lambda\right|}{n}}{2(1+\left|\lambda\right|)}\right]^{\frac{1}{2}}.\label{eq:-3}\end{equation}
}

\begin{flushleft}
\textit{In particular,}
\par\end{flushleft}

\def\theequation{${10}$}\textit{\begin{equation}
\left[\frac{1+r}{2}\left(1-\frac{1}{n}\right)\right]^{\frac{1}{2}}\sqrt{\frac{n}{1-r}}\leq\mathcal{I}_{n,\, r}\left(H^{2},\, B_{2,\,2}^{\frac{1}{2}}\right)\leq\left(1+r+\frac{1}{\sqrt{n}}+\frac{1-r}{n}\right)^{\frac{1}{2}}\sqrt{\frac{n}{1-r}},\label{eq:-3}\end{equation}
}

\def\theequation{${11}$}\begin{equation}
\sqrt{\frac{\frac{1+r}{2}}{1-r}}\leq\liminf_{n\rightarrow\infty}\frac{\mathcal{I}_{n,\, r}\left(H^{2},\, B_{2,\,2}^{\frac{1}{2}}\right)}{\sqrt{n}}\leq\limsup_{n\rightarrow\infty}\frac{\mathcal{I}_{n,\, r}\left(H^{2},\, B_{2,\,2}^{\frac{1}{2}}\right)}{\sqrt{n}}\leq\sqrt{\frac{1+r}{1-r}},\label{eq:-3}\end{equation}
\textit{and }

\def\theequation{${12}$}\begin{equation}
\frac{\sqrt{2}}{2}\leq\liminf_{r\rightarrow1^{-}}\liminf_{n\rightarrow\infty}\sqrt{\frac{1-r}{n}}\mathcal{I}_{n,\, r}\left(H^{2},\, B_{2,\,2}^{\frac{1}{2}}\right)\leq\limsup_{r\rightarrow1^{-}}\limsup_{n\rightarrow\infty}\sqrt{\frac{1-r}{n}}\mathcal{I}_{n,\, r}\left(H^{2},\, B_{2,\,2}^{\frac{1}{2}}\right)\leq\sqrt{2}.\label{eq:-3}\end{equation}

\vspace{0.1cm}

In the next Section, we first give some definitions introducing the
main tools used in the proofs of Theorem A and Theorem B. After that,
we prove these theorems.

\section*{Proofs of Theorems A and B}

From now on, if $\sigma=\left\{ \lambda_{1},\,...,\,\lambda_{n}\right\} \subset\mathbb{D}$
is a finite subset of the unit disc, then \[
B_{\sigma}={\displaystyle \prod_{j=1}^{n}}b_{\lambda_{j}}\]
is the corresponding finite Blaschke product where $b_{\lambda}=\frac{\lambda-z}{1-\overline{\lambda}z},$
$\lambda\in\mathbb{D}$. In Definitions 1, 2, 3 and in Remark 4 below,
$\sigma=\left\{ \lambda_{1},\,...,\,\lambda_{n}\right\} $ is a sequence
in the unit disc $\mathbb{D}$ and $B_{\sigma}$ is the corresponding
Blaschke product.

\begin{flushleft}
\textbf{Definition 1.}\textit{ Malmquist family. }For $k\in[1,\, n]$,
we set $f_{k}=\frac{1}{1-\overline{\lambda_{k}}z},$ and define the
family $\left(e_{k}\right)_{1\leq k\leq n}$, (which is known as Malmquist
basis, see {[}N1, p.117{]}), by
\par\end{flushleft}

\def\theequation{${13}$}\begin{equation}
e_{1}=\frac{f_{1}}{\left\Vert f_{1}\right\Vert _{2}}\,\,\,\mbox{and}\,\,\, e_{k}=\left({\displaystyle \prod_{j=1}^{k-1}}b_{\lambda_{j}}\right)\frac{f_{k}}{\left\Vert f_{k}\right\Vert _{2}}\,,\label{eq:}\end{equation}
for $k\in[2,\, n]$; we have $\left\Vert f_{k}\right\Vert _{2}=\left(1-\vert\lambda_{k}\vert^{2}\right)^{-1/2}.$

\begin{flushleft}
\textbf{Definition 2.}\textit{ The model space $K_{B_{\sigma}}$.
}We define $K_{B_{\sigma}}$ to be the $n$-dimensional space:
\par\end{flushleft}

\def\theequation{${14}$}\begin{equation}
K_{B_{\sigma}}=\left(B_{\sigma}H^{2}\right)^{\perp}=H^{2}\ominus B_{\sigma}H^{2}.\label{eq:}\end{equation}

\begin{flushleft}
\textbf{Definition 3.}\textit{ The orthogonal projection ~$P_{B_{\sigma}}$on
$K_{B_{\sigma}}.$ }We define $P_{B_{\sigma}}$ to be the orthogonal
projection of $H^{2}$ on its $n$-dimensional subspace $K_{B_{\sigma}}.$
\par\end{flushleft}

\begin{flushleft}
\textbf{Remark 4.} The Malmquist family $\left(e_{k}\right)_{1\leq k\leq n}$
corresponding to $\sigma$ is an orthonormal basis of $K_{B_{\sigma}}.$
In particular,
\par\end{flushleft}

\def\theequation{${15}$}\begin{equation}
P_{B_{\sigma}}=\sum_{k=1}^{n}\left(\cdot,\, e_{k}\right)_{H^{2}}e_{k},\label{eq:}\end{equation}

\begin{flushleft}
where $\left(\cdot,\,\cdot\right)_{H^{2}}$ means the scalar product
on $H^{2}$.
\par\end{flushleft}

\begin{flushleft}
\textbf{Proof of Theorem A.}
\par\end{flushleft}

\begin{flushleft}
\textit{Proof of (i).} 1) We fist prove the the right-hand side inequality
of (4). Using both Cauchy-Schwarz inequality and the fact that $\widehat{f'}(k)=(k+1)\widehat{f}(k+1)$
for all $k\geq0,$ we get \[
\left\Vert f'\right\Vert _{L_{a}^{2}}^{2}=\sum_{k\geq0}\frac{\left|\widehat{f'}(k)\right|^{2}}{k+1}=\sum_{k\geq0}\frac{(k+1)^{2}\left|\widehat{f}(k+1)\right|^{2}}{k+1}=\]
\[
=\sum_{k\geq1}k\left|\widehat{f}(k)\right|^{2}\leq\left(\sum_{k\geq1}k^{2}\left|\widehat{f}(k)\right|^{2}\right)^{\frac{1}{2}}\left(\sum_{k\geq1}\left|\widehat{f}(k)\right|^{2}\right)^{\frac{1}{2}}=\]
\[
=\left\Vert f'\right\Vert _{H^{2}}\left\Vert f\right\Vert _{H^{2}}\leq\mathcal{C}_{n,\, r}\left(H^{2},\, H^{2}\right)\left\Vert f\right\Vert _{H^{2}}^{2},\]
and hence, \[
\left\Vert f'\right\Vert _{L_{a}^{2}}\leq\sqrt{\mathcal{C}_{n,\, r}\left(H^{2},\, H^{2}\right)}\left\Vert f\right\Vert _{H^{2}},\]
which means \[
\mathcal{C}_{n,\, r}\left(L_{a}^{2},\, H^{2}\right)\leq\sqrt{\mathcal{C}_{n,\, r}\left(H^{2},\, H^{2}\right)}.\]
Then it remains to use {[}Z2, p.2{]}: \[
\mathcal{C}_{n,\, r}\left(H^{2},\, H^{2}\right)\leq\left(1+r+\frac{1}{\sqrt{n}}\right)\frac{n}{1-r},\]
for all $n\geq1$ and $r\in[0,\,1)$. 
\par\end{flushleft}

\begin{flushleft}
2) The proof of the left-hand side inequality of (4) repeates the
one of {[}Z2, (i){]} (for the left-hand side inequality) excepted
that this time, we replace the Hardy norm $\left\Vert \cdot\right\Vert _{H^{2}}$
by the Bergman one $\left\Vert \cdot\right\Vert _{L_{a}^{2}}$. Indeed,
we use the same test function $e_{n}=\frac{\left(1-r^{2}\right)^{\frac{1}{2}}}{1-rz}b_{r}^{n-1}$
(the $n^{th}$ vector of the Malmquist family associated with the
one-point set $\sigma_{n,\, r}=\underbrace{\{r,\, r,...,\, r\}}_{n}$
see Definition 1) and prove by the same changing of variable $\circ b_{r}$
(in the integral on the unit disc $\mathbb{D}$ which defines the
$L_{a}^{2}-$norm) that \[
\left\Vert e_{n}'\right\Vert _{L_{a}^{2}}^{2}=\frac{n}{1-r}\left(1-\frac{1-r}{n}\right),\]
which gives \[
\mathcal{C}_{n,\, r}\left(L_{a}^{2},\, H^{2}\right)\geq\sqrt{\frac{n}{1-r}}\left(1-\frac{1-r}{n}\right)^{\frac{1}{2}}.\]
Here are the details of the proof. We have $e_{n}\in K_{b_{r}^{n}}$
and $\left\Vert e_{n}\right\Vert _{H^{2}}=1,$ (see {[}N1{]}, Malmquist-Walsh
Lemma, p.116). Moreover,\[
e_{n}'=\frac{r\left(1-r^{2}\right)^{\frac{1}{2}}}{\left(1-rz\right)^{2}}b_{r}^{n-1}+(n-1)\frac{\left(1-r^{2}\right)^{\frac{1}{2}}}{1-rz}b_{r}'b_{r}^{n-2}=\]
\[
=-\frac{r}{\left(1-r^{2}\right)^{\frac{1}{2}}}b_{r}'b_{r}^{n-1}+(n-1)\frac{\left(1-r^{2}\right)^{\frac{1}{2}}}{1-rz}b_{r}'b_{r}^{n-2},\]
since $b_{r}'=\frac{r^{2}-1}{\left(1-rz\right)^{2}}$. Then,\[
e_{n}'=b_{r}'\left[-\frac{r}{\left(1-r^{2}\right)^{\frac{1}{2}}}b_{r}^{n-1}+(n-1)\frac{\left(1-r^{2}\right)^{\frac{1}{2}}}{1-rz}b_{r}^{n-2}\right],\]
 and \[
\left\Vert e_{n}'\right\Vert _{L_{a}^{2}}^{2}=\frac{1}{2\pi}\int_{\mathbb{D}}\left|b_{r}'(w)\right|^{2}\left|-\frac{r}{\left(1-r^{2}\right)^{\frac{1}{2}}}\left(b_{r}(w)\right)^{n-1}+(n-1)\frac{\left(1-r^{2}\right)^{\frac{1}{2}}}{1-rw}\left(b_{r}(w)\right)^{n-2}\right|^{2}\mbox{d}m(w)=\]
\[
=\frac{1}{2\pi}\int_{\mathbb{D}}\left|b_{r}'(w)\right|^{2}\left|\left(b_{r}(w)\right)^{n-2}\right|^{2}\left|-\frac{r}{\left(1-r^{2}\right)^{\frac{1}{2}}}b_{r}(w)+(n-1)\frac{\left(1-r^{2}\right)^{\frac{1}{2}}}{1-rw}\right|\mbox{d}m(w),\]
which gives, using the variables $u=b_{r}(w)$,\[
\left\Vert e_{n}'\right\Vert _{L_{a}^{2}}^{2}=\frac{1}{2\pi}\int_{\mathbb{D}}\left|u^{n-2}\right|^{2}\left|-\frac{r}{\left(1-r^{2}\right)^{\frac{1}{2}}}u+(n-1)\frac{\left(1-r^{2}\right)^{\frac{1}{2}}}{1-rb_{r}(u)}\right|^{2}\mbox{d}m(u).\]
But $1-rb_{r}=\frac{1-rz-r(r-z)}{1-rz}=\frac{1-r^{2}}{1-rz}$ and
$b_{r}'\circ b_{r}=\frac{r^{2}-1}{\left(1-rb_{r}\right)^{2}}=-\frac{\left(1-rz\right)^{2}}{1-r^{2}}$.
This implies \[
\left\Vert e_{n}'\right\Vert _{L_{a}^{2}}^{2}=\frac{1}{2\pi}\int_{\mathbb{D}}\left|u^{n-2}\right|^{2}\left|-\frac{r}{\left(1-r^{2}\right)^{\frac{1}{2}}}u+(n-1)\frac{\left(1-r^{2}\right)^{\frac{1}{2}}}{1-r^{2}}(1-ru)\right|^{2}\mbox{d}m(u)=\]
\[
=\frac{1}{\left(1-r^{2}\right)}\frac{1}{2\pi}\int_{\mathbb{D}}\left|u^{n-2}\right|^{2}\left|\left(-ru+(n-1)(1-ru)\right)\right|^{2}\mbox{d}m(u),\]
which gives\[
\left\Vert e_{n}'\right\Vert _{L_{a}^{2}}=\frac{1}{\left(1-r^{2}\right)^{\frac{1}{2}}}\left\Vert \varphi_{n}\right\Vert _{2},\]
where $\varphi_{n}=z^{n-2}\left(-rz+(n-1)(1-rz)\right).$ Expanding,
we get \[
\varphi_{n}=z^{n-2}\left(-rz+n-1+rz-nrz\right)=\]
\[
=z^{n-2}\left(-nrz+n-1\right)=(n-1)z^{n-2}-nrz^{n-1},\]
and \[
\left\Vert e_{n}'\right\Vert _{L_{a}^{2}}^{2}=\frac{1}{\left(1-r^{2}\right)}\left(\frac{(n-1)^{2}}{n-1}+\frac{n^{2}}{n}r^{2}\right)=\frac{1}{\left(1-r^{2}\right)}\left(n(1+r)-1\right)\]
\[
=\frac{n}{\left(1-r\right)(1+r)}\left((1+r)-\frac{1}{n}\right)=\frac{n}{1-r}\left(1-\frac{1-r}{n}\right)\,,\]
which gives \[
\mathcal{C}_{n,\, r}\left(L_{a}^{2},\, H^{2}\right)\,\geq\sqrt{\frac{n}{1-r}}\left(1-\frac{1-r}{n}\right)^{\frac{1}{2}}.\]

\par\end{flushleft}

\begin{flushleft}

\par\end{flushleft}

\begin{flushleft}
\textit{Proof of (ii).} This is again the same proof as {[}Z2, (ii){]}
(the three steps). More precisely in Step 2, we use the same test
function \[
f=\sum_{k=0}^{s+2}(-1)^{k}e_{n-k},\]
(where $s=\left(s_{n}\right)$ is defined in {[}Z2, p.8{]}), and the
same changing of variable $\circ b_{r}$ in the integral on $\mathbb{D}$.
Here are the details of the proof. 
\par\end{flushleft}

\begin{flushleft}
\textbf{Step 1. }We first prove the right-hand-side inequality:\[
\limsup_{n\rightarrow\infty}\frac{1}{\sqrt{n}}\mathcal{C}_{n,\, r}\left(L_{a}^{2},\, H^{2}\right)\,\leq\sqrt{\frac{1+r}{1-r},}\]
which becomes obvious since \[
\frac{1}{\sqrt{n}}\mathcal{C}_{n,\, r}\left(L_{a}^{2},\, H^{2}\right)\,\leq\frac{1}{\sqrt{n}}\sqrt{\mathcal{C}_{n,\, r}\left(H^{2},\, H^{2}\right)}\,,.\]
and \[
\frac{1}{\sqrt{n}}\sqrt{\mathcal{C}_{n,\, r}\left(H^{2},\, H^{2}\right)}\,\rightarrow\sqrt{\frac{1+r}{1-r}},\]
as $n$ tends to infinity, see {[}Z1{]} p. 2.
\par\end{flushleft}

\begin{flushleft}
\textbf{Step 2.} We now prove the left-hand-side inequality:\[
\liminf_{n\rightarrow\infty}\frac{1}{\sqrt{n}}\mathcal{C}_{n,\, r}\left(L_{a}^{2},\, H^{2}\right)\,\geq\sqrt{\frac{1+r}{1-r}.}\]
More precisely, we show that\[
\liminf_{n\rightarrow\infty}\frac{1}{\sqrt{n}}\left\Vert D\right\Vert _{\left(K_{b_{r}^{n}},\,\left\Vert \cdot\right\Vert _{L_{a}^{2}}\right)\rightarrow H^{2}}\geq\sqrt{\frac{1+r}{1-r}}.\]
Let $f\in K_{b_{r}^{n}}.$ Then, \[
f'=\left(f,\, e_{1}\right)_{H^{2}}\frac{r}{\left(1-rz\right)}e_{1}+\sum_{k=2}^{n}(k-1)\left(f,\, e_{k}\right)_{H^{2}}\frac{b_{r}'}{b_{r}}e_{k}+r\sum_{k=2}^{n}\left(f,\, e_{k}\right)_{H^{2}}\frac{1}{\left(1-rz\right)}e_{k}=\]
\[
=r\sum_{k=1}^{n}\left(f,\, e_{k}\right)_{H^{2}}\frac{1}{\left(1-rz\right)}e_{k}+\frac{1-r^{2}}{(1-rz)(z-r)}\sum_{k=2}^{n}(k-1)\left(f,\, e_{k}\right)_{H^{2}}e_{k}=\]
\[
=\frac{r\left(1-r^{2}\right)^{\frac{1}{2}}}{\left(1-rz\right)^{2}}\sum_{k=1}^{n}\left(f,\, e_{k}\right)_{H^{2}}b_{r}^{k-1}+\frac{\left(1-r^{2}\right)^{\frac{3}{2}}}{(1-rz)^{2}(z-r)}\sum_{k=2}^{n}(k-1)\left(f,\, e_{k}\right)_{H^{2}}b_{r}^{k-1}=\]
\[
=-b_{r}'\left[\frac{r}{\left(1-r^{2}\right)^{\frac{1}{2}}}\sum_{k=1}^{n}\left(f,\, e_{k}\right)_{H^{2}}b_{r}^{k-1}+\frac{\left(1-r^{2}\right)^{\frac{1}{2}}}{z-r}\sum_{k=2}^{n}(k-1)\left(f,\, e_{k}\right)_{H^{2}}b_{r}^{k-1}\right].\]
Now using the change of variables $v=b_{r}(u),$ we get\[
\left\Vert f'\right\Vert _{L_{a}^{2}}^{2}=\int_{\mathbb{D}}\left|b_{r}'(u)\right|^{2}\left|\frac{r}{\left(1-r^{2}\right)^{\frac{1}{2}}}\sum_{k=1}^{n}\left(f,\, e_{k}\right)_{H^{2}}b_{r}^{k-1}+\frac{\left(1-r^{2}\right)^{\frac{1}{2}}}{u-r}\sum_{k=2}^{n}(k-1)\left(f,\, e_{k}\right)_{H^{2}}b_{r}^{k-1}\right|^{2}\mbox{d}u=\]
\[
=\int_{\mathbb{D}}\left|\frac{r}{\left(1-r^{2}\right)^{\frac{1}{2}}}\sum_{k=1}^{n}\left(f,\, e_{k}\right)_{H^{2}}v^{k-1}+\frac{\left(1-r^{2}\right)^{\frac{1}{2}}}{b_{r}(v)-r}\sum_{k=2}^{n}(k-1)\left(f,\, e_{k}\right)_{H^{2}}v^{k-1}\right|^{2}\mbox{d}v.\]
Now, $b_{r}-r=\frac{r-z-r(1-rz)}{1-rz}=\frac{z(r^{2}-1)}{1-rz},$
which gives\[
\left\Vert f'\right\Vert _{L_{a}^{2}}^{2}=\int_{\mathbb{D}}\left|\frac{r}{\left(1-r^{2}\right)^{\frac{1}{2}}}\sum_{k=1}^{n}\left(f,\, e_{k}\right)_{H^{2}}v^{k-1}+\frac{\left(1-r^{2}\right)^{\frac{1}{2}}}{v(r^{2}-1)}(1-rv)\sum_{k=2}^{n}(k-1)\left(f,\, e_{k}\right)_{H^{2}}v^{k-1}\right|^{2}\mbox{d}v=\]
\[
=\frac{1}{1-r^{2}}\int_{\mathbb{D}}\left|r\sum_{k=1}^{n}\left(f,\, e_{k}\right)_{H^{2}}v^{k-1}-(1-rv)\sum_{k=2}^{n}(k-1)\left(f,\, e_{k}\right)_{H^{2}}v^{k-2}\right|^{2}\mbox{d}v=\]
\[
=\frac{1}{1-r^{2}}\int_{\mathbb{D}}\left|r\sum_{k=0}^{n-1}\left(f,\, e_{k+1}\right)_{H^{2}}v^{k}-(1-rv)\sum_{k=0}^{n-2}(k+1)\left(f,\, e_{k+2}\right)_{H^{2}}v^{k}\right|^{2}\mbox{d}v.\]
Thus, 
\par\end{flushleft}

\def\theequation{${16}$}\begin{equation}
\frac{1}{^{\left\Vert f\right\Vert _{H^{2}}\sqrt{n(1+r)}}}\left[\left\Vert (1-rv)\sum_{k=0}^{n-2}(k+1)\left(f,\, e_{k+2}\right)_{H^{2}}v^{k}\right\Vert _{L_{a}^{2}}+\left\Vert r\sum_{k=0}^{n-1}\left(f,\, e_{k+1}\right)_{H^{2}}v^{k}\right\Vert _{L_{a}^{2}}\right]\geq\label{eq:-7-1}\end{equation}

\begin{flushleft}
\[
\geq\sqrt{\frac{1-r}{n}}\frac{\left\Vert f'\right\Vert _{L_{a}^{2}}}{\left\Vert f\right\Vert _{H^{2}}}\geq\]
\[
\geq\frac{1}{^{\left\Vert f\right\Vert _{H^{2}}\sqrt{n(1+r)}}}\left[\left\Vert (1-rv)\sum_{k=0}^{n-2}(k+1)\left(f,\, e_{k+2}\right)_{H^{2}}v^{k}\right\Vert _{L_{a}^{2}}-\left\Vert r\sum_{k=0}^{n-1}\left(f,\, e_{k+1}\right)_{H^{2}}v^{k}\right\Vert _{L_{a}^{2}}\right].\]
Now, \[
(1-rv)\sum_{k=0}^{n-2}(k+1)\left(f,\, e_{k+2}\right)_{H^{2}}v^{k}=\]
\[
=\sum_{k=0}^{n-2}(k+1)\left(f,\, e_{k+2}\right)_{H^{2}}v^{k}-r\sum_{k=0}^{n-2}(k+1)\left(f,\, e_{k+2}\right)_{H^{2}}v^{k+1}=\]
\[
=\sum_{k=0}^{n-2}(k+1)\left(f,\, e_{k+2}\right)_{H^{2}}v^{k}-r\sum_{k=1}^{n-1}k\left(f,\, e_{k+1}\right)_{H^{2}}v^{k}=\]
\[
=\left(f,\, e_{2}\right)_{H^{2}}+2\left(f,\, e_{3}\right)_{H^{2}}v+\sum_{k=2}^{n-2}\left[(k+1)\left(f,\, e_{k+2}\right)_{H^{2}}-rk\left(f,\, e_{k+1}\right)_{H^{2}}\right]v^{k}+\]
\[
-r\left[\left(f,\, e_{2}\right)_{H^{2}}v+(n-1)\left(f,\, e_{n}\right)_{H^{2}}v^{n-1}\right]=\]
\[
=\left(f,\, e_{2}\right)_{H^{2}}+\left[\left(f,\, e_{3}\right)_{H^{2}}-r\left(f,\, e_{2}\right)_{H^{2}}\right]v+\sum_{k=2}^{n-2}\left[(k+1)\left(f,\, e_{k+2}\right)_{H^{2}}-rk\left(f,\, e_{k+1}\right)_{H^{2}}\right]v^{k}+\]
\[
-r(n-1)\left(f,\, e_{n}\right)_{H^{2}}v^{n-1},\]
which gives
\par\end{flushleft}

\def\theequation{${17}$}\begin{equation}
\left\Vert (1-rv)\sum_{k=0}^{n-2}(k+1)\left(f,\, e_{k+2}\right)_{H^{2}}v^{k}\right\Vert _{L_{a}^{2}}^{2}=\label{eq:-7-1-1}\end{equation}

\begin{flushleft}
\[
=\left|\left(f,\, e_{2}\right)_{H^{2}}\right|^{2}+\frac{1}{2}\left|\left(f,\, e_{3}\right)_{H^{2}}-r\left(f,\, e_{2}\right)_{H^{2}}\right|^{2}+\]
\[
+\frac{1}{n}r^{4}(n-1)^{2}\left|\left(f,\, e_{n}\right)_{H^{2}}\right|^{2}+\sum_{k=2}^{n-2}\left|\left(f,\, e_{k+2}\right)_{H^{2}}-\frac{rk}{k+1}\left(f,\, e_{k+1}\right)_{H^{2}}\right|^{2}.\]
On the other hand, 
\par\end{flushleft}

\def\theequation{${18}$}\begin{equation}
\left\Vert r\sum_{k=0}^{n-1}\left(f,\, e_{k+1}\right)_{H^{2}}v^{k}\right\Vert _{L_{a}^{2}}\leq r\left(\sum_{k=0}^{n-1}\frac{1}{k+1}\left|\left(f,\, e_{k+1}\right)_{H^{2}}\right|^{2}\right)^{1/2}\leq r\left\Vert f\right\Vert _{H^{2}},\label{eq:-7}\end{equation}

\begin{flushleft}
Now, let $s=\left(s_{n}\right)$ be a sequence of even integers such
that\[
\mbox{lim}{}_{n\rightarrow\infty}s_{n}=\infty\:\mbox{and}\; s_{n}=o(n)\;\mbox{as}\; n\rightarrow\infty.\]
Then we consider the following function $f$ in $K_{b_{r}^{n}}$:
\[
f=\sum_{k=0}^{s+2}(-1)^{k}e_{n-k}.\]
Applying (17) with such an $f$, we get \[
\left\Vert (1-rv)\sum_{k=0}^{n-2}(k+1)\left(f,\, e_{k+2}\right)_{H^{2}}v^{k}\right\Vert _{L_{a}^{2}}^{2}=\]
\[
=r^{4}\frac{(n-1)^{2}}{n}+\]
\[
+\sum_{l=2}^{n-2}(n-l+1)\left|\left(f,\, e_{n-l+2}\right)_{H^{2}}-\frac{r(n-l)}{n-l+1}\left(f,\, e_{n-l+1}\right)_{H^{2}}\right|^{2},\]
setting the change of index $l=n-k$ in the last sum. This finally
gives\[
\left\Vert (1-rv)\sum_{k=0}^{n-2}(k+1)\left(f,\, e_{k+2}\right)_{H^{2}}v^{k}\right\Vert _{L_{a}^{2}}^{2}=\]
\[
=r^{4}\frac{(n-1)^{2}}{n}+\sum_{l=2}^{s+1}(n-l+1)\left|1+\frac{r(n-l)}{n-l+1}\right|^{2}=\]
\[
=r^{4}\frac{(n-1)^{2}}{n}+\sum_{l=2}^{s+1}(n-l+1)\left[1+r\left(1-\frac{1}{n-l+1}\right)\right]^{2},\]
and\[
\left\Vert (1-rv)\sum_{k=0}^{n-2}(k+1)\left(f,\, e_{k+2}\right)_{H^{2}}v^{k}\right\Vert _{L_{a}^{2}}^{2}\geq\]
\[
\geq r^{4}\frac{(n-1)^{2}}{n}+(s+1-2+1)(n-(s+1)+1)\left[1+r\left(1-\frac{1}{n-(s+1)+1}\right)\right]^{2}=\]
\[
=r^{4}\frac{(n-1)^{2}}{n}+s(n-s)\left[1+r\left(1-\frac{1}{n-s}\right)\right]^{2}.\]
In particular,
\par\end{flushleft}

\[
\left\Vert (1-rv)\sum_{k=0}^{n-2}(k+1)\left(f,\, e_{k+2}\right)_{H^{2}}v^{k}\right\Vert _{L_{a}^{2}}^{2}\geq s(n-s)\left[1+r\left(1-\frac{1}{n-s}\right)\right]^{2}.\]
Now, since $\left\Vert f\right\Vert _{H^{2}}^{2}=s_{n}+3,$ we get\[
\liminf_{n\rightarrow\infty}\frac{1}{n\left\Vert f\right\Vert _{H^{2}}^{2}}\left\Vert (1-rv)\sum_{k=0}^{n-2}(k+1)\left(f,\, e_{k+2}\right)_{H^{2}}v^{k}\right\Vert _{2}^{2}\geq\]
\[
\geq\liminf_{n\rightarrow\infty}\frac{1}{n\left\Vert f\right\Vert _{H^{2}}^{2}}\left\Vert f\right\Vert _{H^{2}}^{2}\left(n-\left\Vert f\right\Vert _{H^{2}}^{2}\right)\left[1+r\left(1-\frac{1}{n-s}\right)\right]^{2}=\]
\[
=\lim_{n\rightarrow\infty}\left(1-\frac{s_{n}}{n}\right)\left[1+r\left(1-\frac{1}{n-s}\right)\right]^{2}=(1+r)^{2}.\]
On the other hand, applying (18) with this $f,$ we obtain\[
\lim_{n\rightarrow\infty}\frac{1}{\sqrt{n}\left\Vert f\right\Vert _{H^{2}}}\left\Vert r\sum_{k=0}^{n-1}\left(f,\, e_{k+1}\right)_{H^{2}}v^{k}\right\Vert _{L_{a}^{2}}=0.\]
Thus, we can conclude passing after to the limit as $n$ tends to
$+\infty$ in (16), that\[
\liminf_{n\rightarrow\infty}\sqrt{\frac{1-r}{n}}\frac{\left\Vert f'\right\Vert _{L_{a}^{2}}}{\left\Vert f\right\Vert _{H^{2}}}=\frac{1}{\sqrt{1+r}}\liminf_{n\rightarrow\infty}\frac{1}{^{\left\Vert f\right\Vert _{H^{2}}\sqrt{n}}}\left\Vert (1-rv)\sum_{k=0}^{n-2}(k+1)\left(f,\, e_{k+2}\right)_{H^{2}}v^{k}\right\Vert _{L_{a}^{2}}\geq\]
\[
\geq\frac{1+r}{\sqrt{1+r}}=\sqrt{1+r},\]
and\[
\liminf_{n\rightarrow\infty}\sqrt{\frac{1-r}{n}}\left\Vert D\right\Vert _{K_{b_{r}^{n}}\rightarrow H^{2}}\geq\liminf_{n\rightarrow\infty}\sqrt{\frac{1-r}{n}}\frac{\left\Vert f'\right\Vert _{L_{a}^{2}}}{\left\Vert f\right\Vert _{H^{2}}}\geq\sqrt{1+r}.\]
\textbf{Step 3. Conclusion.} Using both \textbf{Step 1 }and\textbf{
Step 2}, we get\textbf{ }\[
\limsup_{n\rightarrow\infty}\sqrt{\frac{1-r}{n}}\mathcal{C}_{n,\, r}\left(L_{a}^{2},\, H^{2}\right)=\liminf_{n\rightarrow\infty}\sqrt{\frac{1-r}{n}}\mathcal{C}_{n,\, r}\left(L_{a}^{2},\, H^{2}\right)\,=1+r,\]
which means that the sequence $\left(\frac{1}{\sqrt{n}}\mathcal{C}_{n,\, r}\left(L_{a}^{2},\, H^{2}\right)\right)_{n\geq1}$
is convergent and \[
\lim_{n\rightarrow\infty}\frac{1}{\sqrt{n}}\mathcal{C}_{n,\, r}\left(L_{a}^{2},\, H^{2}\right)=\sqrt{\frac{1+r}{1-r}}.\]

\begin{flushright}
$\square$
\par\end{flushright}

\begin{flushleft}
\textbf{Proof of Theorem B.}
\par\end{flushleft}

\textit{Proofs of inequality (8) and of the right-hand side inequality
of (10).} Let $\sigma$ be a sequence in $\mathbb{D},$ and $B=B_{\sigma}$
the finite Blaschke product corresponding to $\sigma$. If $f\in H^{2},$
we use the same function $g$ as in {[}Z3{]} which satisfies $g_{|\sigma}=f_{|\sigma}.$
More precisely, let $g=P_{B}f\in K_{B}$ (see Definitions 2, 3 and
Remark 4 above for the definitions of $K_{B}$ and $P_{B}$). Then
$g-f\in BH^{2}$ and using the definition of $\mathcal{C}_{n,\, r}\left(L_{a}^{2},\, H^{2}\right),$
\[
\left\Vert g'\right\Vert _{L_{a}^{2}}^{2}\leq\left(\mathcal{C}_{n,\, r}\left(L_{a}^{2},\, H^{2}\right)\right)^{2}\left\Vert g\right\Vert _{H^{2}}^{2}.\]
Now applying the identity (2) to $g$ we get \[
\left\Vert g\right\Vert _{B_{2,\,2}^{\frac{1}{2}}}^{2}\leq\left[\left(\mathcal{C}_{n,\, r}\left(L_{a}^{2},\, H^{2}\right)\right)^{2}+1\right]\left\Vert g\right\Vert _{H^{2}}^{2}.\]
Using the fact that $\left\Vert g\right\Vert _{H^{2}}=\left\Vert P_{B}f\right\Vert _{H^{2}}\leq\left\Vert f\right\Vert _{H^{2}},$
we finally get \[
\left\Vert g\right\Vert _{B_{2,\,2}^{\frac{1}{2}}}\leq\left[\left(\mathcal{C}_{n,\, r}\left(L_{a}^{2},\, H^{2}\right)\right)^{2}+1\right]^{\frac{1}{2}}\left\Vert f\right\Vert _{H^{2}},\]
and as a result,\[
I\left(\sigma,\, H^{2},\, B_{2,\,2}^{\frac{1}{2}}\right)\leq\left[\left(\mathcal{C}_{n,\, r}\left(L_{a}^{2},\, H^{2}\right)\right)^{2}+1\right]^{\frac{1}{2}}.\]
It remains to apply the right-hand side inequality of (4) in Theorem
A to prove the right-hand side one of (10).

\textit{Proof of inequality (9).} 1) We use the same test function\[
f=\sum_{k=0}^{n-1}(1-\vert\lambda\vert^{2})^{\frac{1}{2}}b_{\lambda}^{k}\left(1-\overline{\lambda}z\right)^{-1},\]
as the one used in the proof of {[}Z3, Theorem B{]} (the lower bound,
page 11 of {[}Z3{]}). $f$ being the sum of $n$ elements of $H^{2}$
which are an orthonormal family known as Malmquist's basis (associated
with $\sigma_{n,\,\lambda}=\underbrace{\{\lambda,\lambda,...,\lambda\}}_{n}$,
see Remark 4 above or {[}N1, p.117{]}) , we have $\Vert f\Vert_{H^{2}}^{2}=n$.
\vspace{0.2cm}

2) Since the spaces $H^{2}$ and $B_{2,\,2}^{\frac{1}{2}}$ are rotation
invariant, we have $I\left(\sigma_{n,\,\lambda},H^{2},\, B_{2,\,2}^{\frac{1}{2}}\right)=I\left(\sigma_{n,\,\mu},H^{2},\, B_{2,\,2}^{\frac{1}{2}}\right)$
for every $\lambda,\,\mu$ with $\vert\lambda\vert=\vert\mu\vert=r$.
Let $\lambda=-r$. To get a lower estimate for $\Vert f\Vert_{B_{2,\,2}^{\frac{1}{2}}/b_{\lambda}^{n}B_{2,\,2}^{\frac{1}{2}}}$
consider $g$ such that $f-g\in b_{\lambda}^{n}{\rm Hol}(\mathbb{D})$,
i.e. such that $f\circ b_{\lambda}-g\circ b_{\lambda}\in z^{n}{\rm Hol}(\mathbb{D})$.

\vspace{0.2cm}
 3) First, we notice that\[
\left\Vert g\circ b_{\lambda}\right\Vert _{B_{2,\,2}^{\frac{1}{2}}}^{2}=\left\Vert (g\circ b_{\lambda})^{'}\right\Vert _{L_{a}^{2}}^{2}+\left\Vert g\circ b_{\lambda}\right\Vert _{H^{2}}^{2}=\left\Vert b_{\lambda}.(g'\circ b_{\lambda})\right\Vert _{L_{a}^{2}}^{2}+\left\Vert g\circ b_{\lambda}\right\Vert _{H^{2}}^{2}=\]
\[
=\int_{\mathbb{D}}\left|b_{\lambda}(u)\right|^{2}\left|g'(b_{\lambda}(u))\right|^{2}du+\left\Vert g\circ b_{\lambda}\right\Vert _{H^{2}}^{2}=\int_{\mathbb{D}}\left|g'(w)\right|^{2}dw+\left\Vert g\circ b_{\lambda}\right\Vert _{H^{2}}^{2},\]
 using the changing of variable $w=b_{\lambda}(u)$. We get\[
\left\Vert g\circ b_{\lambda}\right\Vert _{B_{2,\,2}^{\frac{1}{2}}}^{2}=\left\Vert g'\right\Vert _{L_{a}^{2}}^{2}+\left\Vert g\circ b_{\lambda}\right\Vert _{H^{2}}^{2}=\left\Vert g\right\Vert _{B_{2,\,2}^{\frac{1}{2}}}^{2}+\left\Vert g\circ b_{\lambda}\right\Vert _{H^{2}}^{2}-\left\Vert g\right\Vert _{H^{2}}^{2}\,,\]
and\[
\left\Vert g\right\Vert _{B_{2,\,2}^{\frac{1}{2}}}^{2}=\left\Vert g\right\Vert _{H^{2}}^{2}+\left\Vert g\circ b_{\lambda}\right\Vert _{B_{2,\,2}^{\frac{1}{2}}}^{2}-\left\Vert g\circ b_{\lambda}\right\Vert _{H^{2}}^{2}=\]
\[
\geq\left\Vert g\circ b_{\lambda}\right\Vert _{B_{2,\,2}^{\frac{1}{2}}}^{2}-\left\Vert g\circ b_{\lambda}\right\Vert _{H^{2}}^{2}.\]
Now, we notice that

\[
f\circ b_{\lambda}=\sum_{k=0}^{n-1}z^{k}\frac{(1-\vert\lambda\vert^{2})^{\frac{1}{2}}}{1-\overline{\lambda}b_{\lambda}(z)}=\left(1-\vert\lambda\vert^{2}\right)^{-\frac{1}{2}}\left(1+(1-\overline{\lambda})\sum_{k=1}^{n-1}z^{k}-\overline{\lambda}z^{n}\right)=\]

\[
=(1-r^{2})^{-\frac{1}{2}}\left(1+(1+r)\sum_{k=1}^{n-1}z^{k}+rz^{n}\right)\,.\]

\vspace{0.2cm}
 4) Next,\[
\left\Vert g\circ b_{\lambda}\right\Vert _{B_{2,\,2}^{\frac{1}{2}}}^{2}-\left\Vert g\circ b_{\lambda}\right\Vert _{H^{2}}^{2}=\sum_{k\geq1}k\left|\widehat{g\circ b_{\lambda}}(k)\right|^{2}\geq\]
\[
\geq\sum_{k=1}^{n-1}k\left|\widehat{g\circ b_{\lambda}}(k)\right|^{2}=\sum_{k=1}^{n-1}k\left|\widehat{f\circ b_{\lambda}}(k)\right|^{2},\]
since $\widehat{g\circ b_{\lambda}}(k)=\widehat{f\circ b_{\lambda}}(k)\,,\;\forall\, k\in[0,\, n-1].$
This gives\[
\left\Vert g\circ b_{\lambda}\right\Vert _{B_{2,\,2}^{\frac{1}{2}}}^{2}-\left\Vert g\circ b_{\lambda}\right\Vert _{H^{2}}^{2}\geq\frac{1}{1-r^{2}}\left((1+r)^{2}\sum_{k=1}^{n-1}k\right)=\]
\[
=\frac{(1+r)^{2}}{1-r^{2}}\frac{n(n-1)}{2}=\frac{1+r}{1-r}\frac{n(n-1)}{2}=\frac{1+r}{1-r}\frac{(n-1)}{2}\left\Vert f\right\Vert _{H^{2}}^{2},\]
for all $n\geq2$ since $\left\Vert f\right\Vert _{H^{2}}^{2}=n.$
Finally, \[
\left\Vert g\right\Vert _{B_{2,\,2}^{\frac{1}{2}}}^{2}\geq\frac{n}{1-r}\frac{1+r}{2}\left(1-\frac{1}{n}\right)\left\Vert f\right\Vert _{H^{2}}^{2}.\]
In particular, \[
\mathcal{I}_{n,\, r}\left(H^{2},\, B_{2,\,2}^{\frac{1}{2}}\right)\geq\sqrt{\frac{n}{1-r}}\left[\frac{1+r}{2}\left(1-\frac{1}{n}\right)\right]^{\frac{1}{2}}.\]

\begin{flushright}
$\square$
\par\end{flushright}

\subsection*{Some comments}

$\,$

\textbf{a. Extension of Theorem A to spaces $B_{2,\,2}^{s},\, s\geq0$}.
Using the techniques developped in the proof of our Theorem A (combined
with complex interpolation (between Banach spaces) and a reasoning
by induction), it is possible both to precise the sharp numerical
constant $c_{2,\, s}$ in K. Dyakonov's result (3) (mentioned above
in paragraph d. of the Introduction) and to prove the asymptotic sharpness
(at least for $s\in\mathbb{N}\cup\frac{1}{2}\mathbb{N})$ of the right-hand
side inequality of (3). In the same spirit, we would obtain that there
exists a limit: \def\theequation{${19}$}\begin{equation}
\lim_{n\rightarrow\infty}\frac{\mathcal{C}_{n,\, r}\left(B_{2,\,2}^{s-1},\, H^{2}\right)}{n^{s}}=\left(\frac{1+r}{1-r}\right)^{s}.\label{eq:-5-1}\end{equation}
Our Theorem A corresponds to the case $s=\frac{1}{2}.$

\textbf{b. Extension of Theorem B to spaces $B_{2,\,2}^{s},\, s\geq0$.}
The proof of the upper bound in our Theorem B can be extended so as
to give an upper (asymptotic) estimate of the interpolation constant
$\mathcal{I}_{n,\, r}\left(H^{2},\, B_{2,\,2}^{s}\right),\: s\geq0.$
More precisely, applying K. Dyakonov's result (3) (mentioned above
in paragraph d. of the Introduction) we get\def\theequation{${20}$}\begin{equation}
\mathcal{I}_{n,\, r}\left(H^{2},\, B_{2,\,2}^{s}\right)\leq\tilde{c}_{s}\left(\frac{n}{1-r}\right)^{s},\;{\rm with}\;\tilde{c_{s}}\asymp c_{2,\, s},\label{eq:-5-1-1}\end{equation}
where $c_{2,\, s}$ is defined in (3) and precised in (19). Looking
at the above comment 1, $\tilde{c}_{s}\asymp(1+r)^{s}$ for sufficiently
large values of $n$. Our Theorem B corresponds again to the case
$s=\frac{1}{2}.$ In this Theorem B, we prove the sharpness of the
right-hand side inequality in (20) for $s=\frac{1}{2}.$ However,
for the general case $s\geq0,$ the asymptotic sharpness of $\left(\frac{n}{1-r}\right)^{s}$
as $r\rightarrow1^{-}$ and $n\rightarrow\infty$ is less obvious.
Indeed, the key of the proof (for the sharpness) is based on the property
that the Dirichlet norm (the one of $B_{2,\,2}^{1/2}$) is {}``nearly''
invariant composing by an elementary Blaschke factor $b_{\lambda},$
as this is the case for the $H^{\infty}$ norm. A conjecture given
by N. K. Nikolski is the following: \def\theequation{${21}$}\begin{equation}
\mathcal{I}_{n,\, r}\left(H^{2},\, B_{2,\,2}^{s}\right)\asymp\left\{ \begin{array}{c}
\frac{n^{s}}{\sqrt{1-r}}\;\mbox{if}\; s\geq\frac{1}{2}\\
\left(\frac{n}{1-r}\right)^{s}\;\mbox{if}\;0\leq s\leq\frac{1}{2}\end{array},\right.\label{eq:-5-1-1-1}\end{equation}
and is due to the position of the spaces $B_{2,\,2}^{s},\, s\geq0$
with respect to the algebra $H^{\infty}.$

\vspace{0.4cm}

\noun{CMI-LATP, UMR 6632, Université de Provence, 39, rue F.-Joliot-Curie,
13453 Marseille cedex 13, France}

\textit{E-mail address} : rzarouf@cmi.univ-mrs.fr
\end{document}